\documentstyle{article}
\renewcommand{\footnote}[1]{ }
\newcommand{\lcm}[1]{\mbox{\rm lcm}(#1)}

\newcommand{\mod}[1]{\quad \mbox{mod }#1}
\newcommand{\fig}[2]{\begin{figure}[htbp]
        \caption{#1} 
        \begin{center}
                {\small #2}
        \end{center}    
\end{figure}}
\newtheorem{thm}{Theorem}
\newtheorem{prop}[thm]{Proposition}

\newtheorem{lemme}[thm]{Lemma}
\newtheorem{definition}[thm]{Definition}

\newcounter{hour}     
\newcounter{hours}    
\newcounter{minute}   
\newcounter{sixty}    

\newcommand{\setclock} {
\setcounter{hour}{\time}
\setcounter{sixty}{60}
\setcounter{minute}{\time}
\divide \value{hour} by \value{sixty}
\setcounter{hours}{\value{hour}}
\multiply \value{hours} by \value{sixty}
\advance \value{minute} by -\value{hours}
}

\newcommand{\clock}{
\the\value{hour}:\ifnum\the\value{minute}<10 0\fi\the\value{minute}
}

\setclock

\newenvironment{ayo}[9]{\noindent \unitlength=0.53mm
   \special{em:linewidth 0.4pt}
   \linethickness{0.4pt}
   \begin{picture}(71.00,41.00)
 \put(15.00,21.00){\circle{6.00}}
   \put(25.00,21.00){\circle{6.00}}
   \put(35.00,21.00){\circle{6.00}}
   \put(45.00,21.00){\circle{6.00}}
   \put(55.00,21.00){\circle{6.00}}
   \put(65.00,21.00){\circle{6.00}}
   \put(15.00,31.00){\circle{6.00}}
   \put(25.00,31.00){\circle{6.00}}
   \put(35.00,31.00){\circle{6.00}}
   \put(45.00,31.00){\circle{6.00}}
   \put(55.00,31.00){\circle{6.00}}
   \put(65.00,31.00){\circle{6.00}}
   \put(40.00,26.00){\oval(62.00,22.00)[]}
   \put(15.00,12.00){\makebox(0,0)[cc]{\small 7}}
   \put(25.00,12.00){\makebox(0,0)[cc]{\small 6}}
   \put(35.00,12.00){\makebox(0,0)[cc]{\small 5}}
   \put(45.00,12.00){\makebox(0,0)[cc]{\small 4}}
   \put(55.00,12.00){\makebox(0,0)[cc]{\small 3}}
   \put(65.00,12.00){\makebox(0,0)[cc]{\small 2}}
   \put(65.00,40.00){\makebox(0,0)[cc]{\small 1}}
   \put(55.00,40.00){\makebox(0,0)[cc]{\small 0}}
   \put(45.00,40.00){\makebox(0,0)[cc]{\small $-1$}}
   \put(35.00,40.00){\makebox(0,0)[cc]{\small $-2$}}
   \put(25.00,40.00){\makebox(0,0)[cc]{\small $-3$}}
   \put(15.00,40.00){\makebox(0,0)[cc]{\small $-4$}}
   \put(15.00,21.00){\makebox(0,0)[cc]{\small #7}}
   \put(35.00,21.00){\makebox(0,0)[cc]{\small #6}}
   \put(45.00,21.00){\makebox(0,0)[cc]{\small #5}}
   \put(55.00,21.00){\makebox(0,0)[cc]{\small #4}}
   \put(65.00,21.00){\makebox(0,0)[cc]{\small #3}}
   \put(65.00,31.00){\makebox(0,0)[cc]{\small #2}}
   \put(55.00,31.00){\makebox(0,0)[cc]{\small #1}}
    \put(1.00,5.00){\makebox(0,0)[lc]{#8}}
   \put(4.00,26.00){\makebox(0,0)[cc]{#9}}
   }{\end{picture}}
   
\newenvironment{ayononum}[9]{\noindent \unitlength=0.53mm
   \special{em:linewidth 0.4pt}
   \linethickness{0.4pt}
   \begin{picture}(71.00,41.00)
 \put(15.00,21.00){\circle{6.00}}
   \put(25.00,21.00){\circle{6.00}}
   \put(35.00,21.00){\circle{6.00}}
   \put(45.00,21.00){\circle{6.00}}
   \put(55.00,21.00){\circle{6.00}}
   \put(65.00,21.00){\circle{6.00}}
   \put(15.00,31.00){\circle{6.00}}
   \put(25.00,31.00){\circle{6.00}}
   \put(35.00,31.00){\circle{6.00}}
   \put(45.00,31.00){\circle{6.00}}
   \put(55.00,31.00){\circle{6.00}}
   \put(65.00,31.00){\circle{6.00}}
   \put(40.00,26.00){\oval(62.00,22.00)[]}
   \put(15.00,21.00){\makebox(0,0)[cc]{\small #7}}
   \put(35.00,21.00){\makebox(0,0)[cc]{\small #6}}
   \put(45.00,21.00){\makebox(0,0)[cc]{\small #5}}
   \put(55.00,21.00){\makebox(0,0)[cc]{\small #4}}
   \put(65.00,21.00){\makebox(0,0)[cc]{\small #3}}
   \put(65.00,31.00){\makebox(0,0)[cc]{\small #2}}
   \put(55.00,31.00){\makebox(0,0)[cc]{\small #1}}
    \put(1.00,5.00){\makebox(0,0)[lc]{#8}}
   \put(4.00,26.00){\makebox(0,0)[cc]{#9}}
   }{\end{picture}}

\begin{document}


\begin{titlepage}
\title{The Combinatorics of Mancala-Type Games: Ayo, Tchoukaillon, and
$1/\pi$} 
\author{Duane M. Broline\\
Department of Mathematics\\
Eastern Illinois University\\
Charleston, Illinois 61920\\
USA\\
cfdmb@eiu.edu\and
	Daniel E. Loeb\thanks{Author partially supported by URA CNRS 1304,
	EC grant CHRX-CT93-0400, the PRC Maths-Info, and NATO CRG 930554.}\\
LaBRI\\
Universit\'e de Bordeaux I\\
33405 Talence Cedex\\
France\\
loeb@labri.u-bordeaux.fr}
\date{\today, \clock}
\end{titlepage}
\maketitle
\begin{abstract}
Certain endgame considerations in the  two-player Nigerian
Mancala-type game Ayo 
can be identified with the problem of finding winning positions in the
solitaire game Tchoukaillon.
The periodicity of the pit occupancies in
$s$ stone winning positions is determined. Given $n$ pits, the number
of stones in a winning position is found to be asymptotically bounded
by $n^{2}/\pi$. 
\end{abstract}

\section{Introduction}
Around the world, literally hundreds of different kinds of mancala
games have been observed \cite{B2,B3,B5} certain dating back to the Empire
Age of ancient Egypt. Their common features involve cup-shaped
depression called {\em pits} filled with seeds or stones. Players
take turns harvesting 
stones by moving the around the board according to various rules. In
this paper, we will study two variants: The two-player game {Ayo}
played by the Yoruba of western Nigeria, and the solitaire game
Tchoukaillon created by Deledicq and Popova \cite[p.~180]{DP} as a
variant of the game Tchouka played in central Europe. 

We will study certain common unbalanced {Ayo} endgame positions
which we call {\em determined} and show how they are related to
the positions in Tchoukaillon from which a win is possible. For all
$s$, there is a unique such position with $s$ stones. Of course,
certain positions are not realizable on a finite board with a fixed number of 
pits. However, we show that
the number of stones in such a position on a board
with 2$n$ pits (resp.~$n$ in Tchoukaillon) is bounded by approximately
$n^{2}/\pi$. 
Furthermore, we study the actual distribution of stones into pits, and
discover a periodicity  in the contents of the first $k$ pits (with 
respect to the total number of stones) of $\lcm{1,2,\ldots , k+2}$.

\section{Ayo}
The game {\em Ayoyayo} or simply {\em Ayo} is played
on a wooden block 20 inches long, 8 inches wide, and 2 inches
thick. Two rows of six pits each about 3 inches in
diameter are carved in the board. The playing pieces are either dried
palm nuts, or more commonly, the stones of the shrub {\em caesalpina
crista}. The rules \cite{B4} are as follows: 
\begin{itemize}
\item[Set up.] 48 stones are used. Initially, 4 are placed in each of
the 12 pits. (We will generalize the game somewhat allowing 
boards with $2n$ pits and an arbitrary placements of stones.)
\item[Players.] Two players
alternate making moves. Each player's side of the board has $n$ pits.
\item[Objective.] The object of the game is to capture the most stones.
\item[Movement.] To move, a player chooses a non-empty pit from 
his or her side of the board,
and removes all of its stones. The stones are redistributed ({\em sown}),
one per pit, among the pits in a counterclockwise direction
beginning with the pit after the chosen pit.
\item[Odu.] A pit which contains $2n$ or more stones is said to be an
{\em Odu} \cite{B4}. 
If the chosen pit is an {\em Odu}, the redistribution proceeds as usual except that
the initial pit is skipped on each circuit of the board. (None of the
positions we shall consider will contain an {\em Odu}.)
\item[Capture.]If the last pit sown by a player is on the
opponent's side of the board and contains (after having been sown)
two or three stones, then the stones in this pit are captured.  Also
captured are stones in the consecutively preceding pits which meet
these conditions.
\item[End of Game.] At each turn, a player must, if possible, move in such 
a way that his or her opponent has a legal move.  If, on some move, a player
cannot move in such a way to give his or her opponent a legal move, the game is 
over and the player is awarded all remaining stones.  If there are so few
stones on the board that neither player can ever capture, but both players 
will always have a legal move, the game is over and each player is awarded the stones
on his or her own side of the board. For example, if the position is
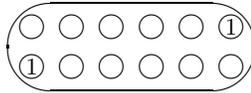
\smallbreak\begin{center}\begin{ayononum}{}{1}{}{}{}{}{1}{}{}
\end{ayononum}\end{center}
no further captures are possible, but each player can always move to give
the opponent a legal move.  In this case, each player is awarded a single stone.
\end{itemize}
The game opens rapidly with both players showing dexterity and skill by
the speed of their movements. However, playing the game well requires
remembering the number of stones in each of the twelve pits, as well
as planning several moves in advance. Thus, the opening game  is both
interesting to watch and difficult to learn.

The endgame is less exciting, but easier to analyze. The latter stages
of the game tend to be dominated by one player. She can
move in such a 
way that her opponent has at all times only one legal move. After this
sequence of moves, only a few stones usually remain, and no further
captures are possible.

We shall analyze a specific type of endgame on a generalized Ayo board with
$2n$ pits.
For reasons which will become apparent when we examine 
Tchoukaillon, the pits will be numbered
clockwise $-n+2,-n+1,\ldots,-1,0,1,\ldots,n,n+1$. (See Figure
\ref{board}.) 
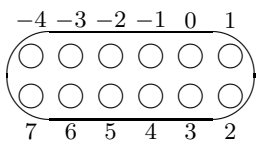
\fig{The standard Ayo board numbering}{
\label{board}\begin{ayo}{}{}{}{}{}{}{}{}{}\end{ayo}}
The two players will be
denoted S (for South) and N (for North).  S makes her plays from 
pits numbered  from $n+1$ down to $2$,
while N makes his plays from  pits numbered from $1$ down to $-n+2$.
Play proceeds from higher 
numbered pits to lower numbered ones (and
 from pit $-n+2$ to pit $n+1$).  

The endgame positions we shall
study are those which satisfy the following definition.
\begin{definition}
A {\em determined position} is an arrangement of stones on a
generalized Ayo board where it is possible for S to move such that
\begin{itemize}
\item S captures at every turn,
\item there is no move from an Odu,
\item \label{waive} after every turn, N has only one stone on his side 
of the board, and
\item all stones are captured by S except one which is awarded to N.
\end{itemize} 
\end{definition}
\fig{\label{F1} The Determined Position with Nine stones}
{\begin{ayo}{ }{1}{2}{2}{4}{ }{ }{N: To move from pit
1}{(a)}\end{ayo}
  \begin{ayo}{1}{ }{2}{2}{4}{ }{ }{S: To capture from pit
2}{(b)}\end{ayo}
  \begin{ayo}{ }{1}{ }{2}{4}{ }{ }{N: To move from pit 1}{(c)}\end{ayo}
\medbreak

  \begin{ayo}{1}{ }{ }{2}{4}{ }{ }{S: To capture from pit 4}{(d)}\end{ayo}
  \begin{ayo}{ }{1}{1}{3}{ }{ }{ }{N: To move from pit 1}{(e)}\end{ayo}
  \begin{ayo}{1}{ }{1}{3}{ }{ }{ }{S: To capture from pit 3}{(f)}\end{ayo}
 \medbreak

  \begin{ayo}{ }{1}{2}{ }{ }{ }{ }{N: To move from pit 1}{(g)}\end{ayo}
  \begin{ayo}{1}{ }{2}{ }{ }{ }{ }{S: To capture from pit 2}{(h)}\end{ayo}
  \begin{ayo}{ }{1}{ }{ }{ }{ }{ }{N: To move from pit 1}{(i)}\end{ayo}
Note:  All captures are to pit 0.  }

Figure \ref{F1} 
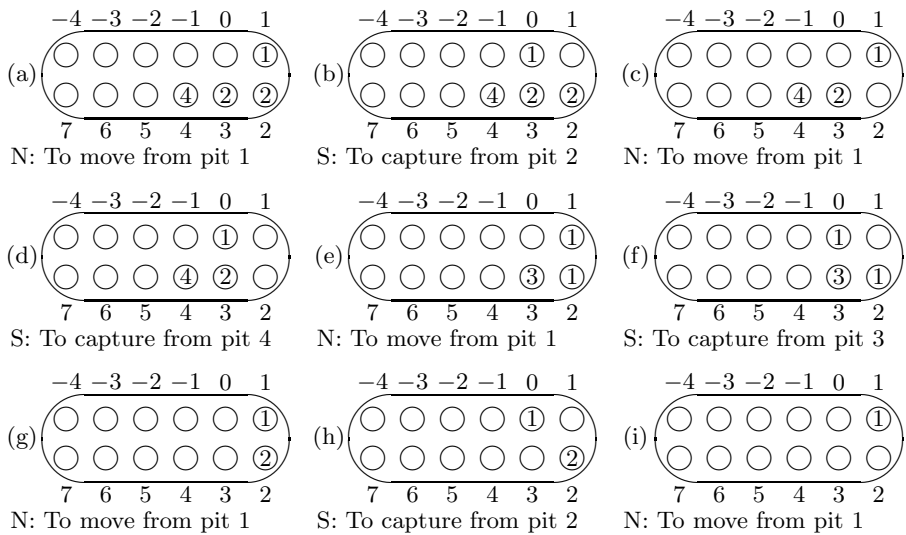
shows a determined
position and the 
subsequent play between the two players on a board with 12 pits.
Initially, there are nine stones and it is N's turn. Eight stones are 
captured by S and one is awarded to N.  It is possible to show 
that a
determined position on a 12 pit {Ayo} board has at most 21
stones.  

It is a simple matter to establish the contents of the pits on N's side
of the board in a determined position.  The study of the contents of the pits
on S's side of the board will be more rewarding.
\begin{lemme}\label{L1}
The stone on N's side of a determined Ayo position must be in pit 1
if N is to move, and in pit 0 if S is to move.
\end{lemme}
{\em Proof:} If S is to move, she must capture and leave only one stone.
Thus, the stone captured must lie in N's second pit (pit 0). Hence, before
N's move, the stone must have been in pit 1.$\Box $

\section{Tchoukaillon}

The game Tchouka is a Russian game  \cite[p.~99]{SL} \cite[p.~42]{SL2}
of possible Paleosiberian or Eskimo origin. It has seemed to
have disappeared several decades ago \cite[p.~99, 180]{DP}. 
The game is played with a number of small pits, dug in the sand, each 
initially containing a certain number of stones and an additional empty 
pit called the {\em Rouma}, {\em Cala} or {\em Roumba}.\footnote{DB: I 
changed a few things here.  Below it is 
important that the number of stones be any number so that the game can be 
played.  Also, I drew a board with an odd number of pits.  It should be 
redrawn if it is important that the number of pits be even!  Also, I have 
changed all occurrences of roumba to Roumba}
(See figure \ref{tchouka}, 
\fig{A Tchouka or Tchoukaillon
board}{\label{tchouka}\begin{picture}(94.00,28.00) 
\put(15.00,21.00){\circle{6.00}}
\put(25.00,21.00){\circle{6.00}}
\put(35.00,21.00){\circle{6.00}}
\put(45.00,21.00){\circle{6.00}}
\put(55.00,21.00){\circle{6.00}}
\put(65.00,21.00){\circle{6.00}}
\put(65.00,8.00){\makebox(0,0)[cc]{2}}
\put(55.00,8.00){\makebox(0,0)[cc]{3}}
\put(45.00,8.00){\makebox(0,0)[cc]{4}}
\put(35.00,8.00){\makebox(0,0)[cc]{5}}
\put(25.00,8.00){\makebox(0,0)[cc]{6}}
\put(15.00,8.00){\makebox(0,0)[cc]{7}}
\put(75.00,21.00){\circle{6.00}}
\put(87.00,21.00){\circle{10.00}}
\put(75.00,8.00){\makebox(0,0)[cc]{1}}
\put(52.00,21.00){\oval(84.00,14.00)[]}
\end{picture}}
although note that sometimes the pits are arranged in a
circle.)  The objective of the game is to put the stones in the Roumba.
As a solitaire game, one may sow the stones from any pit, distributing
them one at a time in the succeeding pits in the direction of the Roumba.
If necessary,
the 
sowing continues
with the pit
opposite
the Roumba.  There are then three possibilities:
\begin{itemize}
\item The last stone falls in the Roumba. The player then continues by
sowing another pit of his choice.
\item The last stone falls in an occupied pit (other than the Roumba).
This pit must then be immediately sown.
\item The last stone falls in an empty pit (other than the Roumba). The
player has lost and the game is over.
\end{itemize}
\footnote{For $n=2$, there is a unique solution. For $n=3$, there is no
solution. For $n=4$, there are 9 solutions. DB: What does solution mean?
It seems non-trivial to prove this result.}
When played as a two-player game, each player continues until he is
forced to place a last stone in an empty pit. At that point, his
opponent moves. Each player attempts to place more stones in the
Roumba than does his adversary.\footnote{DL: This paragraph is new.
Some of the text later is new as well. Maybe I should remove the
comments about two-player Tchouka.

DB: Personally, I find the comments about the two person game to be 
interesting.  This version of a mancala game as certain features that make 
it different from the others.  Most notably the board.}

The solitaire game Tchoukaillon was invented\footnote{We need a reference 
here.} 
\cite[p.~180]{DP} as a variant of Tchouka. In this game, no wrap-around 
moves are allowed,
and the last stone must land in the Roumba.
Thus, pit $i$ may be harvested if and only if 
contains exactly $i$ stones. 

For the purposes of this paper, we will not discuss any alternative
or secondary objectives in the case that this is impossible. We will
thus say that the game is {\em won} if all stones are pocketed. A
position from which a win is possible is said to be a {\em winnable}
position. 

There is a close relationship between Tchoukaillon and Ayo. Suppose 
two people play Ayo on a board of size $2n$ while we focus our
attention on pits  $1,2\dots n+1$ disregarding which player makes a
given move. We would then have the impression that we are watching a
game of Tchoukaillon with pits numbered 1, 2, \dots $n+1$. 
\begin{thm}
There is a one-to-one correspondence between determined positions in
Ayo and winnable positions in Tchoukaillon. To find the Tchoukaillon
position corresponding to a determined Ayo position, ignore pits
$0,-1,-2,\ldots ,-n+2$.
\end{thm}

{\em Proof:} Let $D$ be an Ayo position corresponding to the
Tchoukaillon position $W$. Suppose $D$ is determined. 
\begin{itemize}
\item If pit 1 has a stone,
then by lemma \ref{L1}, player N must move the stone from pit 1 to pit
0. The corresponding move in Tchoukaillon is legal: remove a stone
from pit 1.  
\item If pit 1 is empty, then by lemma \ref{L1}, player S must
capture the stone in pit 0. She does so by harvesting some pit $i$
containing exactly $i$ stones, placing stones in pits 1, 2, \dots
$i-1$, and capturing the stone in pit 0.  The corresponding move in
Tchoukaillon is legal: harvest some pit $i$ containing $i$ stones,
placing stones in pits 1, 2, 3, \dots, $i-1$, and pocketing the remaining
stone in the roumba.
\end{itemize}
Similarly, the legal moves in $W$ correspond to determined moves
in $D$. The objective in Tchoukaillon is to empty the board.
At the end 
of determined Ayo play, the board is empty except for the pit 0.
Thus, $W$ is winning if and only if $D$ is determined.$\Box $  

Given the equivalence between these two notions, we will now
concentrate on Tchoukaillon with the understanding that the results
found will be equally valid for the game of Ayo.

There is a simple strategy by which any feasible win can be forced.
\begin{prop}[\cite{DA,DP}] \label{least} If a win is possible from a
given Tchoukaillon position, the unique winning move must be to
harvest the smallest harvestable pit. 
\end{prop}
{\em Proof:} Suppose pits $i$ and $j$ are both harvestable (i.e., they
contain $i$ and $j$ stones, respectively) and $i<j$. 
If pit $j$ is harvested, then pit $i$ will contain $i+1$ stones. It
would then be ``overfull,'' and could no longer be harvested. Further
play could only increase the occupancy of pit $i.\Box $ 

This strategy can be used \cite{DA,DP} to enumerate a large number of
winning positions.\fig{Winning Positions with up to 24 
stones}{\label{fig1}{
\begin{tabular}{|r|*{8}{c|}|c|} 
\hline
Stones&Pit&Pit&Pit&Pit&Pit&Pit&Pit&Pit&Harvest\\
\multicolumn{1}{|c|}{$s$}&1&2&3&4&5&6&7&8&$h_s$ \\\hline
0& && && && &&{}\\
1&\bf 1&& && && &&1\\ \cline{2-2}
2& &\bf 2&& && && &2\\
3&\bf 1&2&&&&&&&1\\ \cline{2-2}
4& &1&\bf 3&& && &&3\\
5&\bf 1&1&3&& && &&1\\ \cline{2-3}
6& & &2&\bf 4&& &&&4\\
7&\bf 1& &2&4&& &&&1 \\ \cline{2-2}
8& &\bf 2&2&4&& &&&2\\
9&\bf 1&2&2&4&& &&&1\\ \cline{2-2}
10& &1&1&3&\bf 5&&&&5\\
11&\bf 1&1&1&3&5&&&&1\\ \cline{2-4}
12& & & &2&4&\bf 6&&&6\\
13&\bf 1& & &2&4&6&&&1\\ \cline{2-2}
14& &\bf 2& &2&4&6&&&2\\
15&\bf 1&2& &2&4&6&&&1\\ \cline{2-2}
16& &1&\bf 3&2&4&6&&&3\\
17&\bf 1&1&3&2&4&6&&&1\\ \cline{2-3}
18& & &2&1&3&5&\bf 7&&7\\
19&\bf 1& &2&1&3&5&7&&1\\ \cline{2-2}
20& &\bf 2&2&1&3&5&7&&2\\
21&\bf 1&2&2&1&3&5&7&&1\\ \cline{2-2}
22& &1&1&&2&4&6&\bf 8&8\\
23&\bf 1&1&1&&2&4&6&8&1\\ \cline{2-4}
24&&&&\bf 4&2&4&6&8&1\\ \hline
\end{tabular}}}
In fact, the strategy can be applied backwards. That
is to say, given a winning position, one can obtain a winning position
with one more stone in the following manner: Let $i\geq 1$, be the
least number such that pit $i$ is empty. Place $i$ stones in this pit,
and remove one stone from all previous pits. (This can be done since
by definition they are non-empty.) Applying the winning
strategy involves removing the $i$ stones and sowing back 1 stone
into all the pits from which it was removed. We thus obtain an
explicit bijection between winning positions with $s$ stones  and
those with $s+1$. Since there is but one position (winning to be sure)
with no stones, we have the following result.
\begin{thm}
For all $s \geq 0$, there is exactly one winning position involving a
total of $s$ stones.$\Box $
\end{thm}
The winning positions with $s\le 24$ are enumerated in figure
\ref{fig4}. Note in particular that at most 21 stones may appear in a
determined position on a standard Ayo board. Further calculations were
done up to $s=21,286,434$ using a simple SML
program.
\begin{verbatim}
fun ayo carry nil = [carry]
|   ayo carry (0::xs) = carry::xs
|   ayo carry (x::xs) = (x-1)::(ayo (carry+1) xs);
\end{verbatim}
\section{Periodicity}
The columns of Figure \ref{fig1} exhibit a certain periodicity. That is
to say, the number of stones in the first pit depends not  on $s$ but
seemingly on $s$ modulo 2. The contents of the first two pits are periodic 
of period 6 and those of the first three of period 12. 
An extended table\fig{Period of the contents of the first $i$ 
pits}{\label{perfig}\begin{tabular}{r|*{11}{r}} 
$i$&1&2&3&4&5&6&7&8&9&10&11\\ \hline
period&2&6&12&60&60&420&840&2520&2520&27720&27720
\end{tabular}}
suggests surprisingly that the contents of the first four pits and the 
contents of the first five pits have the same periodicity, as do the 
contents of the first eight pits and the contents of the first nine. 
(See Figure \ref{perfig}.)

This periodicity can be established by an analysis of sequences of numbers
determined by the winning positions.  Consider the unique winning position 
with $s$ stones. 
Let $p_{i,s}$ be the number of stones initially in pit $i$.  Clearly,
$p_{i,s}\le i$.  Let $m_{i,s}$ be the number of times pit $i$ must be
harvested in order to win and $b_{i,s}$ be the number of moves of the 
winning strategy which result in a stone being added to pit $i$.  By
convention, $b_{0,s}=s$.
Taking the Roumba to be pit 0, we thus have $b_{0,s}=s$.
Obviously, for $i\geq 1$,
\begin{eqnarray}\label{pimb}
p_{i,s} &=& i m_{i,s} - b_{i,s}.
\end{eqnarray}
There is a one-to-one correspondence between the moves which add a stone to 
pit $i$ and the moves which harvest some pit $j$ for $j > i$.  Hence
$b_{i,s} = \sum_{j>s} m_{j,s}$.  It follows that
\begin{eqnarray}\label{bmb}
b_{i,s} &=& m_{i+1,s} + b_{i+1,s}.
\end{eqnarray}

\begin{prop}\label{periodic}
For all $i$, the sequence of $i$-tuples
\[((p_{1,s};p_{2,s};\cdots;p_{i,s}))_{s\geq 0}\] 
is periodic of period $\lcm{1,2,3,\ldots,i+1}$.
\end{prop}

{\em Proof: } 
We show, by induction on $i$, that not only is the sequence of $i$-tuples 
periodic of period $t=\lcm{1,2,\ldots,i+1}$, but also
$t$ is the smallest positive number such that 
\[ p_{j,t} = 0, \quad j = 1, 2, \ldots, i.\]

The result is trivial for $i=1$.  Assume, by induction, that the result 
holds for all values less than or equal to some $i \ge 1$.
Let $t=\lcm{1,2,\ldots,i+1}$. The inductive hypotheses imply $p_{j,kt}=0$, 
for $j= 1, 2, \ldots, i$ and $k=  1, 2, \ldots.$  Thus, equations 
\ref{pimb} and \ref{bmb} imply
\[jm_{j,kt}=b_{j,kt}=m_{j+1,kt}+b_{j+1,kt}=(j+2)m_{j+1,kt},\quad
j = 1,2,\dots,i-1.\]
Combining these results we get
\[ 2 m_{1,kt} = i(i+1) m_{i,kt}.\]
Since every other move is an addition to pit 1, $2m_{1,kt}=kt.$
Therefore,
$$ \begin{array}{rcll}
p_{i+1,kt}&=&(i+1) m_{i+1,kt} - b_{i+1,kt}\\
&\equiv&(i+1)(m_{i+1,kt}+b_{i+1,kt})&\mod{i+2}\\
&=&(i+1)b_{i,kt}\\
&=& i (i+1)m_{i,kt}\\
&=& 2m_{1,kt}\\
&=& kt 
\end{array}$$

The smallest positive value of $k$ such that 
\[kt\equiv 0\mod{i+2}\quad \] is $k= (i+2)/\gcd(t,i+2)$.  Setting
$q=kt$, we have 
$$\begin{array}{rcccl}
q &=& ((i+2)t)/(\gcd(t,i+2)) &=& \lcm{t,i+2} \\
  &=& \lcm{(\lcm{1, 2, \ldots,i},i+2} &=& \lcm{1, 2, \ldots, i+2}.
\end{array}$$
Thus, we have shown 
\[p_{j,q}=0,\quad i = 1, 2, \ldots i+1\]
and $q$ is the smallest positive multiple of $t$ for which this is true.
By the inductive assumption, we can thus deduce that $q$ is the smallest 
positive integer for which this is true.  In particular, the sequence of 
$i+1$-tuples
\[((p_{1,s};p_{2,s};\cdots;p_{i+1,s}))_{s\ge 0}\]
does not have period less than $q$.

Now, the contents of the first $i+1$ pits in the winning 
position with $q$ stones are the same as the contents of these pits in the 
winning position with no stones. Suppose for some $s$ that the contents of 
the first $i+1$ pits is the same in the unique winning position with $s$ 
stones as in the winning position with $s+q$ stones.  The winning 
positions with $s+1$ and $s+q+1$ stones, respectively, are obtained from 
the corresponding positions with one fewer stones by either adding stones 
to the same pit in both cases or by adding stones to two different pits, 
each having index larger than $i+1$.  In either event, the effect 
upon the first $i+1$ pits is the same and $p_{j,s+1} = p_{j,s+q+1}$, for 
$j = 1, 2, \ldots, i+1.$ We are able to conclude that 
$p_{i+1,s}=p_{i+1,s+q}$ for all $s\ge 0$.  Therefore, the sequence 
of $i+1$-tuples has period $q$. $\Box$

\section{Asymptotics}
There is a winning position for every number of stones given an
unlimited number of pits. However, as in Ayo, if there is a finite
number of pits 
$n$, then not all winning positions are realizable. In particular,
those rows in figure \ref{fig1} of length larger than $n$ can not be
realized. 

Let $s(n)$ denote the smallest number of stones which actually
requires the $n$th pit to win. (Obviously any greater number of stones
will require at least $n$ pits.) 
\begin{figure}[h] \begin{center}
\caption{The minimum number of stones to require $n$ pits is well
approximated by $n^{2}/\pi$}
\label{upto} 
\end{center}
\end{figure}
We will derive the asymptotic formula $s(n)\sim n^{2}/\pi$ from several
interesting observations arising from an examination of the sequences
$p_{i,s}$ and $m_{i,s}$.
\begin{lemme} \label{arithmetic}Given the above notation, 
$p_{i,s}-p_{i-1,s}=(i-1)(m_{i,s}-m_{i-1,s})+2m_{i,s}$.
\end{lemme}
{\em Proof:}
Equations \ref{pimb} and \ref{bmb}.$\Box$

\begin{lemme} \label{noninc}
The sequence $(m_{i,s})_{i=1}^{\infty }$ is non-increasing.
\end{lemme}
{\em Proof:} From Lemma \ref{arithmetic}
\[ p_{i,s} - p_{i-1,s} = (i+1)(m_{i,s}-m_{i-1,s}) + 2m_{i-1,s}.\]
Since $p_{i,s} - p_{i-1,s}\le i$, we have 
$m_{i,s} - m_{i-1,s}\le 0$ and $m_{i,s}\le m_{i-1,s}$, as needed.$\Box$
\begin{thm}\label{mainthm}
As $n$ increases, $s(n)$ grows as $n^{2}/\pi + O(n)$.
\end{thm}
{\em Proof:} 
Let $s$ be fixed.
Define $f(M)$ to be the least $i$ such that $M=m_{i,s}$. 
By lemma \ref{noninc}, $m_{i,s}=M$ if and only if  $i\in I_{M}$ where
$I_{M}=\{f(M),f(M)+1,\ldots,f(M-1)-1 \}$. 
By lemma \ref{arithmetic}, $p_{i,s}-p_{i-1,s}=2M$ for $i,i+1\in
I_{M}$. Thus, the sequence $S_{M}=(p_{i,s})_{i\in I_{M}}$ is a finite
arithmetic sequence with common difference $2M$.
\begin{figure} \begin{center}
\caption{Winning position
pit occupancies $p_{i,1925280}$}
\label{fig4} 
\end{center}
\end{figure}

Now,
\[p_{f(M),s}-p_{f(M)-1,s}=(f(M)-1)(M-m_{f(M)-1,s})+2M.\]
Since $p_{f(M)-1,s}\le f(M)-1$ and $M-m_{f(M)-1,s}\le -1 $,
we see that $p_{f(M),s}$, the leading term of $S_M$, satisfies
\[0 \le p_{f(M),s} \le 2M.\]
 A  similar argument
shows  
$p_{f(M-1)-1,s}$, the final term of $S_{M}$,  satisfies
\[ f(M-1)-2M+1 \le p_{f(M-1)-1,s} \le f(M-1)-1.\]  

To compute the total number of terms $f(M-1)-f(M)$ in $S_{M}$ we
compute the difference of the leading and final terms, divide by the
common difference, and add one.
\begin{eqnarray*}
f(M-1)-f(M) &=& (p_{f(M-1)-1,s}-p_{f(M)})/2M + 1\\
&=& f(M-1)/2M +k_{1}\end{eqnarray*}
where $|k_{1} |\leq 3$.
Hence, 
$$ f(M) = \frac{2M-1}{2M}f(M-1) + k_{2} $$
where $|k_{2} |\leq 3 $,
or explicitly
\begin{eqnarray*}
f(M)&=&\frac{1\times 3\times 5 \times \cdots (2M-1)}{2\times 4 \times
6 \times \cdots \times 2M} n + k_{3}M\\[1mm]
& =& \frac{\frac{1}{2} \frac{3}{2} \ldots
\frac{2M-1}{M}}{M!} n + k_{3}M \\[1mm]
&=& \frac{\Gamma(M+\frac{1}{2})n}{M!\sqrt{\pi}} + k_{3}M
\end{eqnarray*}
where $|k_{3} |\leq 3$,
since $\Gamma(1/2)=\sqrt{\pi}$ taking $n+1=M(0)$.

We must compute $s(n)=\sum _{i=1}^n p_{i,s}.$ Thus, we are led
to the sum of each sequence $I_{M}$. The number of terms has been
already computed to be $f(M-1)/2M+k_{1}$, where $|k_1|\le 3$.
Furthermore, the average term
is 
$$ (p_{f(M-1)-1,s}+p_{f(M),s})/2=f(M-1)/2 + k_{4}M $$
where $|k_{4} |\leq 1$.
Multiplying, we find the sum of $I_{M}$ to be
$f(M-1)^{2}/4M + O(f(M-1))$. We are thus led to calculate
$$ s(n)\sim \sum _{M=1}^{\infty } \frac{\Gamma(M+1/2)^{2}n^{2}}{4\pi
M! (M-1)!} = \frac{n^{2}}{4\pi}\,{}_{2}F_{1}\left({1/2, 1/2 \atop 2
};1\right).$$  
The result follows from Gauss's summation formula \cite{G}, since the
hypergeometric function ${}_{2}F_{1}(1/2,1/2; 2;1)$ is equal to
4.
$\Box $  

\section{Sieves}
By slightly changing the way we think about the game, it is possible to 
generate the numbers $s(n)$ in a different way.  Consider a 
Tchoukaillon board with an infinite number of pits, indexed 1,2,3, 
\dots beginning with the first pit after the {Roumba}.   By the 
remarks following 
Proposition 
\ref{least}, we can play the game
backwards.
The first move is to put one stone in pit 1.  The $k$-th move 
is to first determine the empty pit with lowest index, say pit $j$, 
then add $j$ stones to pit $j$ and remove one stone from each of pits 1, 2, 
\dots, $j-1$.  

Let $h_s$ be the pit to which stones are added on the $s$-th move.
(See Figure \ref{fig1}.)
Clearly $h_s=1$ if and only if $s$ is odd.  Of the remaining terms of 
$(h_s)_{s\ge 1}$ not yet assigned, the first, $h_2$, and every third 
term 
thereafter equals 2.  In general, if all terms with $h_2\le j$ have been 
assigned, the first unassigned term and every $(j+2)$-nd term thereafter 
has the value of $j+1$.  

In the notation of Theorem \ref{mainthm}, $s(n)$ is the first term 
of the sequence $(h_s)_{s\ge 1}$ which is equal to $n$.  The sequence
$(h_{s(n)})_{n\ge 1}$ can be generated by a generalized
``sieve of Eratosthenes'' \cite{D,new}.
(See Figure \ref{seives}.)
\fig{Sieve from Tchoukaillon}{
\label{seives}
\begin{tabbing}
\small\=small\=small\=small\=\kill
    $ L := $ {\tt list of integers, in order, beginning with 1} \\
    $ n := 1$\\
    {\tt repeat forever}\\
          \>$h_n := $ {\tt first unslashed number}\\
          \>{\tt slash out 1st and every $n+1$st unslashed number}\\
          \>$n := n + 1$
\end{tabbing}}
\begin{eqnarray*}
a_n^{(0)} &=& n\\
a_n^{(i+1)} &=& a^{(i)}_{\left\lceil (i+1)n/i \right\rceil}\\
s(i) &=& b_1^{(i)} 
\end{eqnarray*}
for $n,i\ge 0$ where by $\lceil x \rceil$ we denote the smallest
integer greater than or equal to $x$.

Theorem \ref{mainthm} thus proves a conjecture of  Erd\"os and Jabotinsky
\cite[p.~121]{new3}. They had proven that $s(n)=n^2/\pi +O(n^{4/3})$ 
and conjectured that $s(n)=n^2/\pi +O(n)$ based on numerical evidence. 

For comparison, consider the Sieve of Eratosthenes.
\fig{Sieve of Eratosthenes}{\begin{tabbing}
\small\=small\=small\=small\=\kill
    $ L := $ {\tt list of integers, in order, beginning with 2} \\
    $ n := 1$\\
    {\tt repeat forever}\\
          \>$p_n := $ {\tt first unslashed number}\\
          \>{\tt slash out 1st and every $p_{n}$th unslashed number}\\
          \>$n := n + 1$
\end{tabbing}}
It is remarkable 
that the slight difference between these algorithms changes 
the output from a sequence whose $n$th term grows like $n^2/\pi$ to 
one whose $n$th term, by the Prime Number Theorem, grows like $n(\log
n)$.

\section{Remarks}
A Tchoukaillon pit $i$ can hold up to $i$ stones without overfilling.
Thus, a $n$ pit 
position can hold up to a total of $i(i+1)/2$ stones, and the winning
positions have a occupancy rate of $2/\pi + O(1/n) \sim 63.66\% $.

Similarly, since pit $i$ has $i+1$ different possible
occupancies in a winning Tchoukaillon position, there are $(i+1)!$
concievable combinations of pit occupancies for the first $i$ pits.
However, by Proposition 
\ref{periodic}, we see that of these only $\lcm{1,2,3,\ldots,i+1}$
actually occur  in winning Tchoukaillon positions.

We would like to express our surprise and satisfaction that these
easily stated problems turned out to have solutions involving higher
mathematics (even hypergeometric series). In particular, it was
remarkable seeing the constant $\pi$ appear in a combinatorial
problem. 

In this paper, we only studied determined positions in Ayo and
positions in Tchoukaillon from which a total win is possible. It may
be interested to also study strategies designed to maximize the number
of stones captured.

Several researchers have studied chip-firing games played with a
certain number of stones on the nodes of a directed graph \cite{BL}. A
node may be fired if it contains as many stones as its out-degree. One
stone is sent to each of its neighbors. Such games are very
interesting from a mathematical point-of-view, since they have
surprisingly many invariants despite the wealth of choices play
seemingly offers. It might be of interest to relate the theory of
chip-firing games to Mancala games.

The first author studied Ayo while he was teaching at the University
of Ibadan, Nigeria from 1975 to 1978. 
The second author was introduced to Tchoukaillon at the French national
congress of the Association MATh en JEANS \cite{DA} by a group of
junior high 
school students from Coll\`ege l'Ardilli\`ere de
N\'ezant (Saint Brice sous For\^et, France) and Coll\`ege Pierre de
Ronsard (Montmorency, France). 
We thank Paul Campbell for putting us into contact and encouraging us
to write this paper.

\end{document}